\documentclass{article}
\usepackage[latin1]{inputenc}
\usepackage{amsmath}
\usepackage{amssymb}
\usepackage{amsthm}
\usepackage{mathrsfs}
\allowdisplaybreaks

\begin{document}

\theoremstyle{plain} \newtheorem{Thm}{Theorem}
\newtheorem{prop}{Proposition}
\newtheorem{lem}{Lemma}
\newtheorem{step}{Step}
\newtheorem{de}{Definition}
\newtheorem{obs}{Observation}
\newcommand{\kor}{\mathring }
\newcommand{\korz}{\mathring{\mathbb{Z} }_p }
\newcommand{\korq}{\mathring{Q} }
\newcommand{\korx}{\mathring{x} }
\newcommand{\res}{\restriction }
\newcommand{\kalapz }{\hat{\mathbb{Z}}_p }
\theoremstyle{remark} \newtheorem*{pf}{Proof}
\renewcommand\theenumi{(\alph{enumi})}
\renewcommand\labelenumi{\theenumi}
\renewcommand{\qedsymbol}{}
\renewcommand{\qedsymbol}{\ensuremath{\blacksquare}}
\title{The Cayley isomorphism property for groups of order $8p$}
\author{Gábor Somlai\\
Department of Algebra and Number Theory \\
E\"otv\"os Loránd University \\ 1117 Budapest, Pázmány Péter sétány 1/C , Hungary\\
email: zsomlei@gmail.com\thanks{Research supported by the Hungarian
  Scientific Fund (OTKA), grant no. K84233}}
\maketitle

\begin{abstract}
For every prime $p > 3$ we  prove that $Q \times \mathbb{Z}_p$ and $\mathbb{Z}_2 ^3 \times \mathbb{Z}_p $ are DCI-groups. This result completes the description of CI-groups of order $8p$.

\end{abstract}

\section{Introduction}\label{intro}

Let $G$ be a finite group and $S$ a subset of $G$. The Cayley graph
$Cay(G,S)$ is defined by having the vertex set $G$ and $g$ is adjacent
to $h$ if and only if $g h^{-1} \in S$. The set $S$ is called the
connection set of the Cayley graph $Cay(G,S)$. A Cayley graph
$Cay(G,S)$ is undirected if and only if $S=S^{-1}$, where $S^{-1} =
\left\{ \, s^{-1} \in G \mid s \in S \, \right\} $.
Every right multiplication via elements of $G$ is an automorphism of
$Cay(G,S)$, so the automorphism group of every Cayley graph on $G$
contains a regular subgroup isomorphic to $G$. Moreover, this property
characterises the Cayley graphs of $G$.

It is clear that $Cay(G,S) \cong Cay(G,S^{\mu } )$ for every
$\mu \in Aut(G)$. A Cayley graph $Cay(G,S)$ is said to be a
CI-graph if, for each $T \subset G$, the Cayley graphs $Cay(G,S)$ and
$Cay(G,T)$ are isomorphic if and only if there is an automorphism
$\mu$ of $G$ such that $S^{\mu } =T$. Furthermore, a group $G$
is called a DCI-group if every Cayley graph of $G$ is a CI-graph and
it is called a CI-group if every undirected Cayley graph of $G$ is a CI-graph.

It was proved in \cite{4} that $\left< a,z \mid a^p=1, z^8=1,
  z^{-1}az=a^{-1} \right>$ is a CI-group, though not a DCI-group.
Let $G$ be a DCI-group of order $8p$, where $p$ is odd prime.
It can easily be seen that every subgroup of a DCI-group is also a
DCI-group. It follows that the Sylow 2-subgroup of $G$ can only be the
quaternion group $Q$ of order $8$ or $\mathbb{Z}_2^3$.

If $p>8$ or $p=5$, then by Sylow's Theorem the Sylow $p$-subgroup of $G$
is a normal subgroup therefore $G$ is isomorphic to one of the following
groups: $\mathbb{Z}_2^3 \times \mathbb{Z}_p$, $Q \times
\mathbb{Z}_p$, $\mathbb{Z}_2^3 \ltimes \mathbb{Z}_p$ or $Q \ltimes
\mathbb{Z}_p$.
It was proved in \cite{7} in that $\mathbb{Z}_2^3 \times \mathbb{Z}_p$ is a CI-group with respect to ternary relational structures if $p \ge 11$.
Moreover, Dobson and Spiga proved in  \cite{1} that $\mathbb{Z}_2^3 \times
\mathbb{Z}_p$ is a CI-group with respect to color binary relational structures for all primes $p$ and it is a CI-group with respect to
color ternary relational structures if and only if $p \ne 3
\mbox{ and } 7$.

Spiga proved in \cite{2} that $Q  \times \mathbb{Z}_3$ is not a CI-group with respect to colour ternary relational
structures and the non-nilpotent group $Q \rtimes \mathbb{Z}_3$ is not a CI-group.

If $p=7$, then either the Sylow $7$-subgroup is normal, in which case $G$ is as before,
or $G$ has $8$ Sylow $7$-subgroups, when $G \cong \mathbb{Z}_2^{3} \rtimes\mathbb{Z}_7$.
The non-nilpotent groups above are not DCI-groups, see \cite{8}. We show that the other groups are DCI-groups.

\begin{Thm}\label{fotetel}
For every prime $p>3$ the groups $Q \times \mathbb{Z}_p$
and $\mathbb{Z}_2^3 \times \mathbb{Z}_p$ are DCI-groups.
\end{Thm}

Our paper is organized as follows. In section \ref{sec1} we introduce the notation that will be used throughout this paper. In section \ref{sec2} we collect important ideas that we will use in the proof of Theorem \ref{fotetel}. Section \ref{nagy} contains the proof of Theorem \ref{fotetel} for primes $p>8$ and Section \ref{seckicsi} contains the proof of Theorem \ref{fotetel} for $p=5 \mbox{ and } 7$.
\section{Technical details}\label{sec1}
In this section we introduce some notation.
Let $G$ be a group. We use $H \le G$ to denote that $H$ is a subgroup
of $G$ and by $N_{G} (H)$ and $C_{G} (H)$ we denote the normalizer and
the centralizer of $H$ in $G$, respectively.

Let us assume that the group $H$ acts on the set $\Omega$ and let $G$
be an arbitrary group. Then by $G \wr_{\Omega} H$ we denote the wreath
product of $G$ and $H$. Every element $g \in G \wr_{\Omega} H$ can be
uniquely written as $hk$, where $k \in K = \prod_{\omega \in
  \Omega} G_{\omega}$ and $h \in H$. The group $K = \prod_{\omega \in
  \Omega} G_{\omega}$ is
called the base group of $G \wr_{\Omega} H$ and the elements of $K$ can
be treated as functions from $\Omega$ to $G$. If $g \in G \wr_{\Omega} H$ and $g=hk$
we denote $k$ by $(g)_b$. In order to simplify the notation $\Omega$ will be omitted if it is clear from
the definition of $H$ and we will write $G \wr H$.

The symmetric group on the set $\Omega$ will be denoted by $Sym(\Omega)$. Let $G$ be a permutation group on the
set $\Omega $. For a $G$-invariant partition $\mathcal{B}$ of the set
$\Omega$ we use $G^{\mathcal{B}}$ to denote the permutation group on $\mathcal{B}$ induced by the action of $G$ and
similarly, for every $g \in G$ we denote by $g^{\mathcal{B}}$ the
action of $g$ on the partition $\mathcal{B}$.

For a group $G$, let $\hat{G}$ denote the subgroup of the
symmetric group $Sym(G)$ formed by the elements of $G$ acting by
right multiplication on $G$. For every Cayley graph $\Gamma =Cay(G,S)$ the
subgroup $\hat{G}$ of $Sym(G)$ is contained in $Aut(\Gamma)$.

\begin{de}\label{zart}
Let $G \le Sym(\Omega)$ be a permutation group. Let
\begin{equation*}
\begin{split}
G^{(2)} = \left\{ \pi \in Sym(\Omega) \bigg| \begin{matrix}
\forall a,b \in
  \Omega \mbox{  } \exists g_{a,b} \in G \mbox{ with}  &\pi
  (a)=g_{a,b} (a)
 \mbox{ and } \\ &\pi (b)=g_{a,b} (b) \end{matrix} \right\} \mbox{.}
\end{split}
\end{equation*}
We say that $G^{(2)}$ is the $2$-closure of the permutation group $G$.
\end{de}

\begin{lem}\label{2zart}
Let $\Gamma$ be a graph. If $G \le Aut(\Gamma )$, then
$G^{(2)} \le Aut(\Gamma )$.
\end{lem}

\section{Basic ideas}\label{sec2}
In this section we collect some results and some important ideas that
we will use in the proof of Theorem \ref{fotetel}.

We begin with a fundamental lemma that we will use all along this paper.
\begin{lem}[Babai \cite{5}]\label{babai}
$Cay(G,S)$ is a CI-graph if and only if for every regular subgroup
$\mathring{G}$ of $Aut(Cay(G,S) )$ isomorphic to $G$ there is a
$\mu \in Aut(Cay(G,S))$ such that ${\mathring{G}}^{\mu} = \hat{G}$.
\end{lem}

We introduce the following definition.
\begin{de}
\begin{enumerate}
\item
We say that a Cayley graph $Cay(G,S)$ is a CI$^{(2)}$-graph iff for every regular subgroup $\mathring{G}$ of
$Aut(Cay(G,S) )$ isomorphic to $G$ there is a $\sigma \in \langle \mathring{G}, \hat{G} \rangle ^{(2)}$
such that ${\mathring{G}}^{\sigma} = \hat{G}$.
\item A group $G$ is called a DCI$^{(2)}$-group if for every $S \subset G$ the Cayley graph $Cay(G,S)$ is a CI$^{(2)}$-graph.
\end{enumerate}
\end{de}

Let us assume that $A= Aut(Cay(G,S) ) \le Sym(8p)$ contains two copies of
regular subgroups, $\mathring{Q} \times \mathring{\mathbb{Z}}_p$
and $\hat{Q} \times \hat{\mathbb{Z}}_p$.
By Sylow's theorem we may assume that $\mathring{\mathbb{Z}}_p$ and
$\hat{\mathbb{Z}}_p$ are in the same Sylow $p$-subgroup $P$ of
$Sym(8p)$. If $p>8$, then $P$ is isomorphic
to $\mathbb{Z}_p^8$. Moreover, $P$ is generated by $8$ disjoint $p$-cycles. It follows
that both $\mathring{Q}$ and $\hat{Q}$ normalize $P$ so we may assume
that $\mathring{Q}$ and $\hat{Q}$ lie in the same Sylow $2$-subgroup of
$N_A (P)$.
Let $P_2$ denote a Sylow $2$-subgroup of $Sym(8)$. It
is also well known that $P_2$ is isomorphic to the automorphism group
of the following graph $\Delta$:

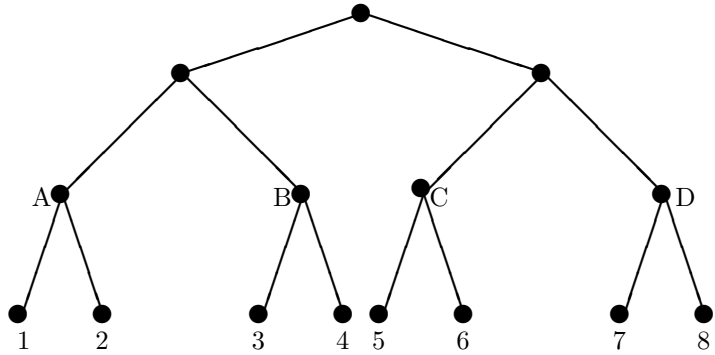
\unitlength 0.8mm
\begin{picture}(129,70)
\thicklines
\put(65,60){ \line(-3,-1){30} }
\put(65,60){ \line(3,-1){30} }

\put(35,50){ \line(-1,-1){20} }
\put(35,50){ \line(1,-1){20} }
\put(95,50){ \line(1,-1){20} }
\put(95,50){ \line(-1,-1){20} }

\put(15,30){ \line(-1,-3){7} }
\put(15,30){ \line(1,-3){7} }
\put(55,30){ \line(1,-3){7} }
\put(55,30){ \line(-1,-3){7} }
\put(75,30){ \line(1,-3){7} }
\put(75,30){ \line(-1,-3){7} }
\put(115,30){ \line(1,-3){7} }
\put(115,30){ \line(-1,-3){7} }

\put(66,60){\circle*{3} }
\put(36,50){\circle*{3} }
\put(96,50){\circle*{3} }

\put(16,30){\circle*{3} }
\put(13,31){\makebox(0,0)[t]{A}}
\put(56,30){\circle*{3} }
\put(53,31){\makebox(0,0)[t]{B}}
\put(76,31){\circle*{3} }
\put(79,31){\makebox(0,0)[t]{C}}
\put(116,30){\circle*{3} }
\put(120,31){\makebox(0,0)[t]{D}}

\put(9,10){\circle*{3} }
\put(23,10){\circle*{3} }

\put(49,10){\circle*{3} }
\put(63,10){\circle*{3} }
\put(69,10){\circle*{3} }
\put(83,10){\circle*{3} }
\put(109,10){\circle*{3} }
\put(123,10){\circle*{3} }

\put(10,7){\makebox(0,0)[t]{1}}
\put(23,7){\makebox(0,0)[t]{2}}
\put(49,7){\makebox(0,0)[t]{3}}
\put(63,7){\makebox(0,0)[t]{4}}
\put(69,7){\makebox(0,0)[t]{5}}
\put(83,7){\makebox(0,0)[t]{6}}
\put(109,7){\makebox(0,0)[t]{7}}
\put(123,7){\makebox(0,0)[t]{8}}

\end{picture}

Every automorphism of $\Delta$ permutes the leaves of the graph and the
permutation of the leaves determines the automorphism, therefore
$Aut(\Delta )$ can naturally be embedded into $Sym(8)$.

It is easy to see that the same holds if we change $Q \times
\mathbb{Z}_p$ to $\mathbb{Z}_2^3 \times \mathbb{Z}_p$.

\begin{lem}\label{leszamlalas}
\begin{enumerate}
\item
There are exactly two regular subgroups of $P_2$ which are
isomorphic to $Q$.
\item
There are exactly two regular subgroups of $P_2$ which are
isomorphic to $\mathbb{Z}_2^3$.

\end{enumerate}

\end{lem}
\begin{pf}
\begin{enumerate}
\item
Let $Q$ be a regular subgroup of $Aut(\Delta)$ isomorphic to the quaternion group with
generators $i$ and $j$.
For every $1 \le m \le 4$ there is a $q_m \in Q$ such that $q_m
(2m-1)=2m$. These are automorphisms of $\Delta$ so $q_m (2m)=2m-1$ and
hence the order of $q_m$ is $2$. There is only one involution in $Q$
so $q_m =i^2$ for every $1 \le m \le 4$ and this fact determines
completely the action of $i^2$ on $\Delta$.

We can assume that $i(1) =3$. Such an isomorphism of $\Delta$ fixes
setwise $\left\{ 1,2,3,4 \right\}$ so we have that $i(3)=2$, $i(2)=4$
and $i(4)=1$ since $i$ is of order $4$. Using again the fact that $Q$
is regular on $\Delta$ and $i^2(5)=6$, we get that there are two
choices for the action of $i$: $i =(1324)(5768)$ or $i =(1324)(5867)$.

We can also assume that $j(1)=5$. This implies that $j(5)=j^2 (1)=i^2(1)=2$,
and $j(2)=6$ since $j \in Aut(\Delta)$ and $j(6)=1$. The action of $i$
determines the action of $j$ on $\Delta$ since $iji=j$. Applying this
to the leaf $3$ we get that $j(3)=8$ if $i =(1324)(5768)$ and $j(3)=7$
if $i =(1324)(5867)$ so there is no more choice for the action of
$j$. Finally, $i$ and $j$ generate $Q$ and this gives the result.
\item
Let us assume that $x \in \mathbb{Z}_2^3$ such that $x(1)=2$
A fixed point free automorphism of $\Gamma$ of order $2$ which maps $1$ to $2$ will map
$3$ to $4$.  There is an $y \in\mathbb{Z}_2^3$ such that $y(1) =5$.
Such an automorphism of $\Gamma$ maps $2$ to $6$ so  we have
that $x(5) =6$ since $x$ and $y$ commute. This determines $x$
completely so we have that $x$ =(12)(34)(56)(78).

We have two possibilities for $y(3)$. If $y(3) =7$, then
$y=(15)(26)(37)(48)$ and if $y(3) =8$, then $y=(15)(26)(38)(47)$.
The third generator of the group $\mathbb{Z}_2^3$ which maps $1$ to
$3$ is determined by $x$ and $y$ since $\mathbb{Z}_2^3$ is abelian.

\end{enumerate}
\qed
\end{pf}

The previous proof also gives the following.
\begin{lem}\label{megad}
\begin{enumerate}
\item\label{megada}
 The following two pairs of permutations generate the two regular
      subgroups of $Aut(\Delta) \le Sym(8)$ isomorphic to $Q$:
\[ i_1 =(1324)(5768) \mbox{, } j_1 =(1526)(3748) \] and
\[ i_2 =(1324)(5867) \mbox{, } j_2 =(1526)(3847)\]
\item\label{megadb} The elements of these regular subgroups of $Aut(\Delta)$ are the following:
\[ \begin{matrix}
Q_l: & Q_r: \\
id & id\\
(12)(34)(56)(78) & (12)(34)(56)(78)\\
(1324)(5768) & (1324)(5867) \\
(1423)(5867) & (1423)(5768) \\
(1526)(3748) & (1526)(3847) \\
(1625)(3847) & (1625)(3748) \\
(1728)(3546) & (1728)(3645) \\
(1827)(3645) & (1827)(3546) \\
\end{matrix} \]
Using the following identification $Q_l$ and $Q_r$ act  on $Q$ by left-multiplication and right-multiplication, respectively:
\[ \begin{matrix}
1 & 2 & 3 & 4 & 5 & 6 & 7 & 8 \\
1 & -1 & i & -i & j & -j & k & -k \\
\end{matrix} \mbox{.}\]
\item The following permutations generate two regular subgroups of
      $Aut(\Delta ) \le Sym(8)$ isomorphic to $\mathbb{Z}_2^3$.
\\ $A_1$ is generated by:
\[x_1 =(12)(34)(56)(78) \mbox{, } x_2=(13)(24)(57)(68) \mbox{, } x_3
=(15)(26)(37)(48) \] and $A_2$ is generated by:
\[y_1 =(12)(34)(56)(78) \mbox{, } y_2=(13)(24)(58)(67) \mbox{, } y_3
=(15)(26)(38)(47) \mbox{.} \]

\end{enumerate}

\end{lem}

\begin{lem}\label{25}
Let us assume that $G_1 \le P_2$ is generated by two different regular subgroups $Q_a$ and $Q_b$ of $Aut(\Delta)$
which are isomorphic to $Q$ and
$G_2 \le P_2$ is generated by two different regular subgroups $A_1$ and $A_2$ of $Aut(\Delta)$ which are isomorphic to $\mathbb{Z}_2^{3}$. Then $G_1 = G_2$.
\end{lem}
\begin{pf}
It is clear that $|P_2| = |Aut(\Delta)| = 2^{7}$. One can see using Lemma \ref{megad} \ref{megada} that $G_1$ and $G_2$ are generated by even permutations.
Both $G_1$ and $G_2$ induce an action on the set $V=\left\{A,~ B,~ C,~ D \right\}$ which is a set of vertices of $\Delta$ and
it is easy to verify that every permutation of $V$ induced by $G_1$ and $G_2$ is even.
This shows that $G_1$ and $G_2$ are contained in a subgroup of $P_2$ of cardinality $2^{5}$.

Lemma \ref{megad} \ref{megadb} shows that $|Q_a \cap Q_b| =2$ and one can also check that $|A_1 \cap A_2| =2$.
This gives $|G_1| \ge 2^{5}$ and $|G_2| \ge 2^{5}$, finishing the proof of Lemma \ref{25}.
\qed
\end{pf}

\begin{prop}\label{ujprop}
\begin{enumerate}
\item The quaternion group $Q$ is a $DCI^{(2)}$-group.
\item $\mathbb{Z}_2^3$ is a $DCI^{(2)}$-group.
\end{enumerate}
\end{prop}
\begin{pf}
\begin{enumerate}
\item Let $Q_a$ and $Q_b$ be two regular subgroups of $Sym(8)$ isomorphic to the quaternion group $Q$.
By Sylow's theorem we may assume that $Q_a$ and $Q_b$ lie in the same Sylow $2$-subgroup of $H=\langle Q_a, Q_b \rangle$. Since every Sylow $2$-subgroup of $H$ is contained in a Sylow $2$-subgroup of $Sym(8)$, we may assume that $Q_a$ and $Q_b$ are subgroups of $Aut(\Delta)$.

 Our aim is to find an element $\pi \in \langle Q_a, Q_b \rangle^{(2)}$ such that $Q_a ^{\pi} = Q_b$ so let us assume that $Q_a \ne Q_b$.
Using Lemma \ref{megad} \ref{megada} we may also assume that $Q_a$ and $Q_b$ are generated by the permutations
$(1324)(5768), (1526)(3748)$ and $(1324)(5867),(1526)(3847)$, respectively.
Lemma \ref{megad} \ref{megadb} shows that $H$ contains the following three permutations:
\begin{equation*} \begin{split}
(12)(34) &=(1324)(5768)(1324)(5867) \\ (12)(56) &=(1526)(3748)(1526)(3847) \\ (12)(78) &=(1728)(3546)(1728)(3645) \mbox{.}
\end{split} \end{equation*}

Now one can easily see that the permutation $(12)$ is in $H^{(2)}$. Finally, it is also easy to check using Lemma \ref{megad} \ref{megadb} that $Q_a^{(12)} = Q_b$.

\item Let $A_1$ and $A_2$ be two regular subgroups of $Sym(8)$ isomorphic to $\mathbb{Z}_2^3$. Let $H'$ denote the group generated by $A_1$ and $A_2$.
    Similarly to the previous case we may assume that $A_1$ and $A_2$ are different regular subgroups of $Aut(\Delta)$. By Lemma \ref{megad} $A_1$ and $A_2$ are generated by the permutations $x_1 =(12)(34)(56)(78)$, $x_2= (13)(24)(57)(68)$, $x_3 = (15)(26)(37)(48)$ and $x_1 = y_1= (12)(34)(56)(78)$, $y_2 = (13)(24)(58)(67)$, $y_3= (15)(26)(38)(47)$, respectively.

    By Lemma \ref{25} the group $H'$ contains the permutations $(12)(34)$, $(12)(56)$ and $(12)(78)$. Therefore $H'$ contains the permutation $(12)$ which conjugates $A_1$ to $A_2$ since $(12)$ centralizes $x_1$ and we also have $(12)x_2 (12) = y_2 y_1$ and $(12) x_3 (12)= y_1 y_3$, finishing the proof of Proposition \ref{ujprop}. \qed
\end{enumerate}
\end{pf}

\begin{de}
Let $\Gamma$ be an arbitrary graph and $A, B \subset V(
\Gamma)$ such that $A \cap B = \emptyset$. We write $A \sim B$ if one of the following four possibilities holds:
\begin{enumerate}
\item For every $a \in A$ and $b \in B$ there is an edge from $a$ to
      $b$ but there is no edge from $b$ to $a$.
\item For every $a \in A$ and $b \in B$ there is an edge from $b$ to
      $a$ but there is no edge from $a$ to $b$.
\item For every $a \in A$ and $b \in B$ the vertices $a$ and $b$ are
      connected with an undirected edge.
\item There is no edge between $A$ and $B$.
\end{enumerate}
We also write $A \nsim B$ if none of the previous four possibilities holds.
\end{de}

\begin{lem}\label{koszoru}
Let $A$, $B$ be two disjoint subsets of cardinality $p$ of a
graph. We write $A \cup B = \mathbb{Z}_p \cup \mathbb{Z}_p $. Let
us assume that $\hat{\mathbb{Z}}_p$ acts naturally on $A \cup B$
and for a generator $\mathring{a} $ of the cyclic group
$\mathring{Z}_p$ the action of $\mathring{a}$ is defined by
$\mathring{a}(a_1,a_2)=(a_1 +b, a_2 +c)$ for some $b, c \in
\mathbb{Z}_p$.
\begin{enumerate}
\item\label{kosza} If $b=c$, then the action of $\mathbb{Z}_p$ and
      $\mathring{Z}_p$ on $A \cup B$ are the same.
\item\label{koszb} If $A \nsim B$, then $b=c$.
\item\label{koszc} If $A \sim B$, then every $\pi \in Sym(A \cup B)$
      which fixes $A$ and $B$ setwise is an automorphism of the graph
      defined on $A \cup B$ if $\pi \restriction A \in Aut(A)$ and
      $\pi \restriction B \in Aut(B)$.
\end{enumerate}
\end{lem}
\begin{pf}
These statements are obvious.
\qed
\end{pf}

\section{Main result for $p > 8$}\label{nagy}

In this section we will prove that $Q \times \mathbb{Z}_p$ and
$\mathbb{Z}_2^3 \times \mathbb{Z}_p$ are DCI-groups if $p>8$.
We will first prove it for $Q \times \mathbb{Z}_p$ and then we
will repeat the argument for the case of $\mathbb{Z}_2^3 \times
\mathbb{Z}_p$.

\begin{prop}\label{pr1}
For every prime $p>8$, the group $Q \times \mathbb{Z}_p$ is a
DCI-group.
\end{prop}
Our technique is based on Lemma \ref{babai} so we have to fix a
Cayley graph $\Gamma =Cay(Q \times \mathbb{Z}_p,S)$. Let
$A=Aut(\Gamma) $ and $ \mathring{G}=\mathring{Q} \times
\mathring{\mathbb{Z} }_p $ be a regular subgroup of $ A $ isomorphic
to $Q \times \mathbb{Z}_p$. In order to prove Proposition
\ref{pr1} we  have to find an $\alpha \in A$ such that
$\mathring{G}^{\alpha} = \hat{G} =\hat{Q} \times \hat{\mathbb{Z} }_p $
what we will achieve in three steps.

\subsection{Step 1}\label{step1}

We may assume $\hat{\mathbb{Z}}_p$ and
$\mathring{\mathbb{Z}}_p$
lie in the same Sylow p-subgroup $P$ of $Sym(8p)$. Then both $\mathring{Q}$ and $\hat{Q}$ are subgroups of $N_{Sym(8p) } (P) \cap A$ so we may assume that $\mathring{Q}$ and $\hat{Q}$ lie in the same Sylow $2$-subgroup of $N_{Sym(8p) } (P) \cap A$ which is contained in a Sylow $2$-subgroup of $A$.

The Sylow $p$-subgroup $P$ gives a partition $\mathcal{B}=\left\{ B_1, B_2, \ldots ,B_8
\right\}$ of the vertices of $\Gamma$, where $\left| B_i \right| =p$
for every $i=1, \ldots ,8$ and $\mathcal{B}$ is $P$-invariant.
It is easy to see that $\mathcal{B}$ is invariant under the action of
$\hat{Q}$ and $\mathring{Q}$ and hence $\langle \hat{G} , \mathring{G}
\rangle \le Sym(p) \wr Sym(8)$.
Moreover, both  $\mathring{G}$ and
$\hat{G}$ are regular so $\mathring{Q}$ and $\hat{Q}$ induce regular
action on $\mathcal{B}$ which we denote by $Q_1$ and $Q_2$, respectively. The assumption that $\mathring{Q}$ and $\hat{Q}$ lie in the same Sylow $2$-subgroup of $A$
 implies that $Q_1$ and $Q_2$ are in the same Sylow $2$-subgroup of $Sym(8)$.

\subsection{Step 2}\label{step2}
Let us assume that $Q_1 \ne Q_2$. We intend to find an element $\alpha \in A$ such that $\left( \mathring{Q}^{\alpha} \right) ^{\mathcal{B} }  = Q_2$.

Using Lemma \ref{megad}\ref{megadb} we can assume that $\korq$ is
 generated by the permutations $\mathring{i}$ and $\mathring{j}$ such that $\mathring{i}$ and $\mathring{j}$ induce the permutations
$(B_1 B_3 B_2 B_4)(B_5 B_7 B_6 B_8)$ and $(B_1B_5B_2B_6)(B_3B_7B_4B_8)$, respectively.
Similarly, $\hat{Q}$ is generated by $\hat{i}$ and $\hat{j}$ with $\hat{i}^{\mathcal{B} } = (B_1 B_3 B_2 B_4)(B_5 B_8 B_6 B_7)$
 and $\hat{j}^{\mathcal{B} } =(B_1B_5B_2B_6)(B_3B_8B_4B_7)$.

We define a graph $\Gamma _0$ on $\mathcal{B}$ such that $B_i$ is
connected to $B_j$ if and only if $B_i \nsim B_j$. This is an undirected
graph with vertex set $\mathcal{B}$ and both $Q_1$ and $Q_2$ are regular subgroups of $Aut(\Gamma_0)$.
It follows that $\Gamma_0$ is a Cayley graph of the quaternion group of order $8$.

\begin{de}
\begin{enumerate}
\item For a pair $(B_i, B_j) \in \mathcal{B}^{2}$ we write $B_i \equiv B_j$ if either there exists a path
$C_1, C_2, \ldots ,C_n$ in $\Gamma_0$ such that $C_1 =B_1$, $C_n =B_2$ or $i=j$.
\item For a pair $(B_i, B_j) \in \mathcal{B}^{2}$ we write $B_i \not\equiv B_j$ if $B_i \equiv B_j$ does not hold.
\item If both $H$ and $K$ are subsets of the vertices of $\Gamma_0$ such that $H \cap K = \emptyset$ and
for every $B_i \in H$, $B_j \in K$ we have $B_i \not\equiv B_j$, then we write $H \not\equiv K$.
\end{enumerate}
\end{de}

\begin{obs}
\begin{enumerate}
\item
The relation $\equiv$ defines an equivalence relation on $\mathcal{B}$. The equivalence classes defined by the relation $\equiv$ will be called equivalence classes.
\item Since $Q_1$ acts transitively on $\mathcal{B}$
we have that the size of the equivalence classes defined by the relation $\equiv$ divides $8$.
\end{enumerate}
\end{obs}

We can also define a colored graph $\Gamma_1$ on $\mathcal{B}$ by coloring the edges of the complete directed graph on $8$ points. $B_i$ is connected to $B_j$ with the same color as $B_i'$ is connected to $B_j'$ in $\Gamma_1$ if and only if there exists a graph isomorphism $\phi$ from $B_i \cup B_j$ to $B_i' \cup B_j'$ such that $\phi(B_i) = B_i'$ and $\phi(B_j) = B_j'$. The graph $\Gamma_1$ is a colored Cayley graph of the quaternion group. Moreover, both $Q_1$ and $Q_2$ act regularly on $\Gamma_1$. Using the fact that $Q$ has property $DCI^{(2)}$ it is clear that there exists an $\alpha' \in \langle Q_1, Q_2 \rangle ^{(2)} \le Aut(\Gamma_1)$ such that $Q_2^{\alpha'} =Q_1$.
We would like to lift $\alpha'$ to an automorphism $\alpha$ of $\Gamma$ such that $\alpha^{\mathcal{B} } = \alpha'$.

\begin{enumerate}
\item
Let us assume first that $\Gamma_0$ is a connected graph.
\begin{lem}\label{osszefuggo}
\begin{enumerate}
\item\label{het} $\mathring{Q} \times \korz \le \hat{\mathbb{Z}}_p \wr Sym(8)$.
\item\label{e} If $\mathring{Q} \times \korz \le \hat{\mathbb{Z}}_p \wr Sym(8)$, then for every $\mathring{q} \in \mathring{Q}$ we have $(\mathring{q} ) _b =id$.
\end{enumerate}
\end{lem}
\begin{pf}
\begin{enumerate}
\item We first prove that $\hat{\mathbb{Z}}_p = \mathring{\mathbb{Z}}_p$.
Let $x$ and $y$ generate $\hat{\mathbb{Z} }_p$ and $\korz$, respectively. We can assume that $x \res{B_1} =y \res{B_1}$.
Using Lemma \ref{koszoru}\ref{koszb} we get that $x \res{B_i} =y \res{B_i}$ if there exists a path in $\Gamma_0$ from $B_1$ to $B_i$.
This shows that $x=y$ since $\Gamma_0$ is connected. Moreover, $\korq \times \korz \le \hat{\mathbb{Z}}_p \wr Sym(8) $
since the elements of $\korz$
and the elements of $\korq$ commute.

\item
Let $A' = A \cap \hat{\mathbb{Z}}_p \wr Sym(8)$. We have already assumed that $\korq$ and $\hat{Q}$ lie in the same Sylow $2$-subgroup of $A'$. Let $\mathring {q}$ be an arbitrary element of $\korq$. For every $(a,u) \in Q \times \mathbb{Z}_p$ we have $\mathring{q}(a,u) = (b, u+t)$ for some $b \in Q$ and $t \in \mathbb{Z}_p$, where $t$ only depends on $\mathring{q}$ and $a$ since $\mathring{q} \le \hat{\mathbb{Z}}_p \wr Sym(8)$. The permutation group $\hat{G}$ is transitive, hence there exist $\hat{q}_1, \hat{q}_2 \in \hat{Q}$ such that $\hat{q}_1(1,u)=(a,u)$ and $\hat{q}_2(b,u+t)=(1,u+t)$. The order of $\hat{q}_2 \mathring{q} \hat{q}_1$
is a power of $2$ since $\hat{q}_2, \mathring{q}, \hat{q}_1$ lie in a Sylow $2$-subgroup. Therefore
 $t =0$ and hence $(\mathring{q})_b = id$.
 \end{enumerate}
 \qed
\end{pf}

Lemma \ref{osszefuggo} says that if $\Gamma_0$ is connected, then $\langle \mathring{Q}, \hat{Q} \rangle \le \hat{Z}_p \wr Sym(8)$ and $(q)_b=id$ for every $q \in \langle \mathring{Q}, \hat{Q} \rangle$.
Therefore we can define $\alpha= \alpha' id_{\mathcal{B} }$ to be an element of the wreath product  $\hat{Z}_p \wr Sym(8)$ and clearly $\alpha' id_{\mathcal{B} }$ is an element of $A$ with $\alpha^{\mathcal{B} } = \alpha'$.

\item
Let us assume that $\Gamma_0$ is the empty graph.

Then Lemma \ref{koszoru}\ref{koszc} shows that every permutation in $\langle Q_1, Q_2 \rangle^{(2)}$ lifts to an automorphism of $\Gamma$.
\item
Let us assume that $\Gamma_0$ is neither connected nor the empty graph.

\begin{obs}\label{involuciok}
If $Q_1 \ne Q_2$, then $ \langle \korq, \hat{Q} \rangle \le A$ contains $\beta_1$, $\beta_2$, $\beta_3$
such that \[ \beta_1^{\mathcal{B} }= (B_1 B_2)(B_3B_4) , ~\beta_2^{\mathcal{B} }= (B_1 B_2)(B_5B_6), ~\beta_3^{\mathcal{B} }= (B_1 B_2)(B_7B_8) \mbox{.}\]
\end{obs}
\begin{pf}
By Lemma \ref{megad} the elements $\beta_1$, $\beta_2$, $\beta_3$ can be generated as products of an element of $\korq$ and $\hat{Q}$.
\qed
\end{pf}
\begin{lem}\label{8}
We claim that $B_{2k-1} \equiv B_{2k}$ for ${k=1,2,3,4}$.
\end{lem}
\begin{pf}
Since $\Gamma_0$ is a Cayley graph and $Q_1$ is transitive on the pairs of the form $(B_{2k-1}, B_{2k})$ it is enough to prove that $B_1 \equiv B_2$.
If $B_1 \nsim B_2$, then $B_1 \equiv B_2$ so we can assume that $B_1 \sim B_2$. Since $\Gamma_0$ is not the empty graph $B_1$ is connected to $B_l$ for some $l>2$. By Observation (\ref{involuciok}) there exists $\beta \in A$ such that $\beta(B_1) =B_2$ and $\beta(B_l) =B_l$. This shows that $B_2 \nsim B_l$ and hence $B_1 \equiv B_2$. \qed
\end{pf}

$\Gamma_0$ is not connected so the size of the equivalence classes defined by $\equiv$ cannot be bigger than $4$. Since $B_1 \equiv B_2$ there exists a partition $H_1 \cup H_2 = \mathcal{B}$ such that $|H_1| =|H_2|=4$, $B_1, B_2 \in H_1$ and $H_1 \not\equiv H_2$.

\begin{lem}\label{9}
There exists $\alpha \in A$ such that $\alpha^{\mathcal{B} } = \alpha'$.
\end{lem}
\begin{pf}
Let us assume first that $H_1= \left\{ B_1, B_2, B_3, B_4 \right\}$. Then we define $\alpha_1$ to be equal to $\beta_2$ on $H_1$ and the identity on $H_2$. Using Lemma \ref{koszoru}\ref{koszb} we get that $\alpha_1$ is in $\langle \korq, \hat{Q} \rangle ^{(2)}$.

If $H_1= \left\{ B_1, B_2, B_5, B_6 \right\}$ or $H_1= \left\{ B_1, B_2, B_7, B_8 \right\}$, then we define $\alpha_2$ by $\alpha_2 \restriction H_1 =\beta_1$ and
$\alpha_2 \restriction H_2 = id$. Lemma \ref{koszoru}\ref{koszb} shows again that $\alpha_2 \in A$.

It is easy to see that $\alpha_1 ^{\mathcal{B}} = \alpha_2 ^{\mathcal{B}} = (B_1 B_2)$. Therefore $A$ contains an element $\alpha$ such that $Q_1^{\alpha^{\mathcal{B} } } =Q_2$.
\qed
\end{pf}

\end{enumerate}

We conclude that we can assume that $Q_1 =Q_2$.

\subsection{Step 3}\label{step3}
Let us now assume that $Q_1 = Q_2$. We intend to find $\gamma \in A$ such that $\mathring{Q} ^{\gamma} = \hat{Q}$.

Let $\hat{x}$ and $\kor{x}$ denote the generators of $\kalapz$ and $\mathring{\mathbb{Z}}_p$, respectively. We may assume that $\hat{x}\restriction B_1 =\kor{x}\restriction B_1$ .
\begin{lem}\label{10}
There exists $\gamma \in A$ such that $\kor{x}^{\gamma} = \hat{x}$.
\end{lem}
\begin{pf}
Let us assume first that $\Gamma_0$ is connected. In this case there is only one equivalence class of size $8$.
It is clear by Lemma \ref{koszoru} \ref{koszb} that $\kor{x} = \hat{x}$.

Let us assume that $\Gamma_0$ is not connected.
In this case there are at least two equivalence classes which we denote by $\mathscr{C}_1, \ldots ,\mathscr{C}_n$. The permutations $\hat{x}$ and $\kor{x}$ are elements of the base group of $\hat{\mathbb{Z}}_p \wr Sym(8)$ and hence they can be considered as functions on $\mathcal{B}$.
By Lemma \ref{koszoru} \ref{koszb} $\kor{x}$ is constant on
every equivalence class and we may assume that $\hat{x}(q,u) =(q, u+1)$ for every $(q,u) \in Q \times \mathbb{Z}_p$. We may also assume that $B_1 \in \mathscr{C}_1$.

For every $1 \le m \le n$ there exists $\mathring{q} _m \in \korq$ such that $\mathring{q}_m ( \mathscr{C} _1 )= \mathscr{C} _m$ and for every $\mathring{q}_m \in \mathring{Q}$ there exists $\hat{q}_m \in \hat{Q}$
such that $\mathring{q}_m^{\mathcal{B}} = \hat{q}_m^{\mathcal{B}}$.
Let $\gamma$ be defined as follows:
\begin{equation*}
\begin{split}
\gamma \restriction \cup \mathscr{C} _1 &= id\\
\gamma \restriction \cup \mathscr{C} _m  &= \kor{q}_m \hat{q}_m^{-1} \mbox{ } \mbox{ for $2 \le m \le n$.}
\end{split}
\end{equation*}

Let $(b,v) \in \mathring{q}_m (B_e)$ with $B_e \in \mathscr{C}_1$ and we denote $\mathring{q}_m^{-1}(b,v)$ by $(a,u)$.
Since $\korx$ is constant on $\mathscr{C}_m$ we have $\korx^{s} (b,v) =(b, v + c_m s)$ for some $c_m$ which only depends on $\mathscr{C}_m$.
Thus $\mathring{q}_m (a, u+s )= (b,v+c_m s  )$ since $\korx$ and $\mathring{q}_m$ commute and $\korx \restriction B_e = \hat{x} \restriction B_e$.
Therefore for every $ w \in \mathbb{Z}_p$ we have
\[ \gamma(b,w) = \mathring{q}_m (a, w) = \mathring{q}_m (a, u+(w-u)) =(b,v+ c_m (w-u) ) \]
for every $(b,w) \in \mathring{q}_m (B_e)$. It is easy to verify that $\gamma^{-1}(b,w) =(b, \frac{w-v +uc_m}{c_m})$
 for every $w \in \mathbb{Z}_p$ which gives
\[ \gamma^{-1} \korx \gamma (b,w) = \gamma^{-1} \korx (b, w c_m +v -u c_m) =\gamma^{-1}
(b, w c_m +v -u c_m +c_m) = (b,w+1) \mbox{.}  \]
It remains to show that $\gamma \in A$.
Let $y$ and $z$ be two points of $Q \times \mathbb{Z}_p$.

If $y$ and $z$ are in the same equivalence class
$\mathscr{C} _m$, then either $\gamma$ is defined on $y$ and
$z$ by $\kor{q}_m \hat{q}_m^{-1}$ which is the element of
the group $\langle \mathring{G}, \hat{G} \rangle \le A$ or
$\gamma(y)=y$ and  $\gamma(z)=z$.

We denote by $B_y$ and $B_z$ the elements of $\mathcal{B}$
containing $y$ and $z$, respectively.
If $y$ and $z$ are not in the same equivalence class,
then $B_y \sim B_z$. The definition of $\gamma$ shows
that $\gamma^{\mathcal{B} }  = id $.
Using Lemma \ref{koszoru} \ref{koszc}
we get that $\gamma \restriction B_y \cup B_z$ is an
automorphism of the induced subgraph of $\Gamma$
on the set $B_y \cup B_z$, which proves that $\gamma \in A$,
finishing the proof of Lemma \ref{10}.
\qed
\end{pf}
Using Lemma \ref{10} we may assume that $\korx = \hat{x}$.
Since $\korx$ and $\mathring{q}$ commute we have $\korq \times \korz \le \hat{\mathbb{Z}}_p \wr Sym(8) $.
Now we can apply Lemma \ref{osszefuggo} which gives $(\kor{q})_b= id$ for every $\kor{q} \in \mathring{Q}$.
This proves that $\korq = \hat{Q}$ since $Q_1 =Q_2$. Therefore $\kor{G} = \hat{G}$, finishing the proof of Proposition \ref{pr1}.

Our method also gives the analogous result for $\mathbb{Z}_2 ^3 \times \mathbb{Z}_p$, what also follows from the theorem of Dobson and Spiga \cite{1}.

\begin{prop}\label{pr2}
For every prime $p>8$, the group $\mathbb{Z}_2^3 \times \mathbb{Z}_p$ is a
DCI-group.
\end{prop}

In order to prove Proposition \ref{pr2} we will mofidy the proof of Proposition \ref{pr1}.
Let $\Gamma$ be a Cayley graph of $G=\mathbb{Z}_2^3 \times \mathbb{Z}_p$ and let $A = Aut (\Gamma)$.
Let $\mathring{G} =\mathring{\mathbb{Z}}_2^3 \times \mathring{\mathbb{Z}}_p$ be a regular subgroup of $A$ isomorphic to $\mathbb{Z}_2^3 \times \mathbb{Z}_p$.
It is enough to prove that there exists $\alpha \in A$ such that $\mathring{G}^{\alpha}=(\mathring{\mathbb{Z}}_2^3 \times \mathring{\mathbb{Z}}_p )^{\alpha}= \hat{\mathbb{Z}}_2^3 \times \hat{\mathbb{Z}}_p =\hat{G}$.

It is easy to verify that the argument of the first step in subsection \ref{step1} only uses the fact that $p>8$.
Therefore there exists a $P$-invariant partition $\mathcal{D}=\left\{ D_1, D_2, \ldots , D_8 \right\}$, where $P$ is a Sylow $p$-subgroup of $Sym(8p)$ containing $\hat{\mathbb{Z}}_p$ and $\mathring{\mathbb{Z}}_p$. We denote be $A_1$ and $A_2$ the regular action on $\mathcal{D}$ induced by
$\hat{\mathbb{Z}}_2^3$ and $\mathring{\mathbb{Z}}_2^3$, respectively.

Let us assume that $A_1 \ne A_2$. We will repeat the argument of Step 2.
Similarly to the definition of $\Gamma_1$ one can define a colored graph
$\Gamma_1'$ on $\mathcal{D}$. Since $\mathbb{Z}_2^3$ is also a $DCI^{(2)}$-group there exists $\beta' \in Aut(\Gamma_1')$ such that $A_2^{\beta'} =A_1$.

One can also define the graph $\Gamma_0'$ using the relation $\equiv$ and similarly to Lemma \ref{osszefuggo} one can prove that if
$\Gamma_0'$ is connected, then there exists $\beta \in A$ such that $\beta^{\mathcal{D} }= \beta'$.

If $\Gamma_0 '$ is the empty graph, then every automorphism of $\Gamma_1'$ lifts to an automorphism of $\Gamma$.

Similarly to Observation \ref{involuciok} the automorphism group $A$ contains
$\delta_1$, $\delta_2$, $\delta_3$
such that \[ \delta_1^{\mathcal{D} }= (D_1 D_2)(D_3 D_4) , ~\delta_2^{\mathcal{D} }= (D_1 D_2)(D_5 D_6), ~\delta_3^{\mathcal{D} }= (D_1 D_2)(D_7 D_8) \mbox{.}\]
since $\langle A_1, A_2 \rangle = \langle Q_1, Q_2 \rangle$ by Lemma \ref{25}.

It is straightforward to check that Lemma \ref{8} and Lemma \ref{9} only uses the existence of the involutions $\beta_1, \beta_2, \beta_3$ so the argument can be repeated using $\delta_1, \delta_2$ and $\delta_3$.
Therefore we may assume that $A_1 = A_2$.

Finally, the proof of Lemma \ref{10} can also be repeated for $\mathbb{Z}_2^3 \times \mathbb{Z}_p$ which gives that the generators of $\hat{\mathbb{Z}}_p$ and $\mathring{\mathbb{Z}}_p$ coincide. Since $A_1 = A_2$ we have $\hat{G}= \mathring{G}$, finishing the proof of Proposition \ref{pr2}. \qed

It is straigtforward to check that the proof of Proposition
\ref{pr1} and Proposition \ref{pr2} only uses the fact that $p>8$ in the first step of the argument. We can formulate this
fact in Proposition \ref{particio}.

\begin{prop}\label{particio}
Let $\Gamma$ be a Cayley graph of $G =Q \times \mathbb{Z}_p$ or $G
=\mathbb{Z}_2^3 \times \mathbb{Z}_p$, where $p$ is an odd
prime and let $\mathring{G} =\mathring{Q} \times \mathring{\mathbb{Z}}_p$ or $\mathring{G} =\mathring{\mathbb{Z}}_2^3 \times \mathring{\mathbb{Z}}_p$ be a regular subgroup of $Aut(\Gamma)$ isomorphic to $G$. Let us assume that there exists a $\langle \hat{G} , \mathring{G} \rangle $-invariant partition $\mathcal{B}=\left\{ B_1, B_2,
  \ldots  ,B_8 \right\}$ of $V(\Gamma)$, where $\left| B_i
\right| =p $ for every $i= \left\{ 1, \dots , 8 \right\}$. In
addition, we assume that $\mathring{\mathbb{Z}}_p$ is a subgroup of the base group of $\hat{\mathbb{Z}}_p wr Sym(\mathcal{B})$. Then
there is an automorphism $\alpha$ of the graph $\Gamma$ such that
$\hat{G}^{\alpha} =\mathring{G}$.
\end{prop}

\section{Main result for $p=5$ and $7$}\label{seckicsi}
In this section we will prove that $Q \times
\mathbb{Z}_5$, $Q \times \mathbb{Z}_7$, $\mathbb{Z}_2^3 \times
\mathbb{Z}_5$ and $\mathbb{Z}_2^3 \times \mathbb{Z}_7$ are CI-groups.

The whole section is based on the paper \cite{4}, so we will only
modify the proof of Lemma 5.4 of \cite{4}.
\begin{prop}\label{kicsi}
Every Cayley graph of $Q \times
\mathbb{Z}_5$, $Q \times \mathbb{Z}_7$, $\mathbb{Z}_2^3 \times
\mathbb{Z}_5$ and $\mathbb{Z}_2^3 \times \mathbb{Z}_7$ is a CI-graph.
\end{prop}
We denote by $R$ one of the groups $Q$ and $\mathbb{Z}_2^3$.
Let $\Gamma$ be a Cayley graph of one of these groups, $A=Aut(
\Gamma )$ and $P$ a Sylow $p$-subgroup of $A$ for $p=5,7$,
respectively. Let us assume that $A$ contains  two copies
of regular subgroups which we denote by $\hat{G} = \hat{R} \times \kalapz$ and $\mathring{G} = \mathring{R} \times \mathring{\mathbb{Z}} _p$.
We can assume that $\Gamma$ is neither the empty nor the complete graph and both
$\hat{\mathbb{Z}}_p$ and $\mathring{\mathbb{Z}}_p$ are contained
in $P$.

It was proved in \cite{4}  that the action
of $A$ on the points of graph $\Gamma$ cannot be primitive so there is
a nontrivial $A$-invariant partition $\mathcal{B}=\left\{ B_0, B_1,
  \ldots , B_{t-1} \right\}$ of $V( \Gamma )=G$. The elements of the
partition $\mathcal{B}$ have the same cardinality since the action of
$A$ is transitive on $\mathcal{B}$ so $\left| B_i \right| <p^2$ for
every $i=0,1, \ldots ,t-1$. The partition $\mathcal{B}$ is $P$-invariant so $P$ acts
on $\mathcal{B}$. Since $P$ is a $p$-group, the length of every orbit of $P$ is a power of $p$.

If the length of every orbit of $P$ on $V(\Gamma)$ is $p$, then it is
clear from Proposition \ref{particio} that $\Gamma$ is a
CI-graph. Therefore $P$ has an orbit $\Lambda \subset G$ such that
$| \Lambda |= p^2$ since $p^3 > |G | $ and the remaining orbits of $P$
have length $p$ since $2p^2> 8p$.

Let $\mathscr{C} = \left\{ C_0, C_1, \ldots , C_{s-1} \right\}$ be an
orbit of $P$ on $\mathcal{B}$ such that $\Lambda \subseteq
\cup_{i=0}^{s-1} C_i$. We may assume that $B_i = C_i$ for $i=0,1, \ldots , s-1$. It is clear that $s$ is a power of $p$. If
$s \ge p^2$, then $\left| \cup_{i=0}^{s-1} C_i \right| \ge 2p^2 >8p$
which is a contradiction.  It follows
that $1<s<p^2$ which implies $s=p$.

For every $i<s$ and every $x \in P$ the following eqalities hold  for
some $j<s$ \[ {(B_i \cap \Lambda)}^x = {B_i}^x \cap {\Lambda}^x
=B_j \cap \Lambda \mbox{.} \] This implies that \[ \left| B_0 \cap
\Lambda \right| = \left| B_i \cap \Lambda \right| \] for every $0 \le i<s$. Therefore
\[ p^2 = \left| \Lambda \right| = \left| \cup_{i=0}^{s-1} \left( B_i \cap
  \Lambda \right) \right| =s \left| B_0 \cap \Lambda \right|  = p \left| B_0 \cap
  \Lambda \right| \mbox{.} \]
This gives $\left| B_0 \cap \Lambda \right| =p$ so $\left| B_0
\right|$ can only be $p$ or $8$ since $\left| B_0 \right| t =8p$ and
both $\left| B_0 \right|$ and $t$ are at least $p$.

If $\left| B_0 \right| = p$, then $\Lambda$ is the union of $p$
elements of the $A$-invariant partition $\mathcal{B}$ and every
orbit $\Lambda '$ of $P$ is an element of the partition
$\mathcal{B}$ if $\Lambda' \ne \Lambda$.
For every orbit $\Lambda ' \ne \Lambda$ of $P$ and for every
$y \in \hat{ \mathbb{Z} }_p \cup \mathring{\mathbb{Z} }_p$ we have $y(\Lambda') = \Lambda'$.
By Proposition \ref{particio} we may assume that there exists an element
$x'$ in $\hat{\mathbb{Z} }_p \cup \mathring{\mathbb{Z} }_p$
such that $x' (B_0) \ne B_0$ and
clearly $x' (B_7) = B_7$ for every $x' \in \hat{ \mathbb{Z} }_p \cup \mathring{\mathbb{Z} }_p$.
Since both $\hat{G}$ and $\mathring{G} $ are regular there exists $a \in C_A(x')$ such that $a(B_0) =B_7$,
which contradicts the fact that $\mathcal{B}$ is $A$-invariant and $B_7$ is an orbit of $P$.

Let us assume that $\left| B_0 \right| =8$ and let $\hat{x}$ and
$\mathring{x}$ generate $\hat{\mathbb{Z} }_p$ and $\mathring{\mathbb{Z} }_p$, respectively.
Since $\hat{G}$ and $\mathring{G}$ are regular we have that neither $\hat{x}^{\mathcal{B}}$
nor $\mathring{x} ^{\mathcal{B} }$ is the identity, while for every $r \in \hat{R} \cup \mathring{R}$
we have $r^{\mathcal{B} } = id$.
Since $\hat{x}$ and $\mathring{x}$ are in the same Sylow
$p$-subgroup of $P$ we may assume that $\hat{x}(B_i) =\mathring{x}
(B_i) = B_{i+1}$
for $i=0,1, \ldots ,p-1$, where the indices are taken modulo $p$.
By Proposition \ref{particio} we may also assume that $\hat{x} \ne \mathring{x}$.

For every $m$ there exists an $l$ such that the action of
$\hat{x}^l \mathring{x} ^{-l}$ is nontrivial on $B_m$ since
$\hat{x} \ne
\mathring{x}$. Therefore $A_{B_m} \restriction B_m$ contains
a regular subgroup and a cycle of length $p$ such that $p >\frac{\left| B_0 \right|}{2}$. A theorem of Jordan says
that such a permutation group is $2$-transitive and hence the induced
subgraph by $B_m$ of $\Gamma$ is the complete or the empty graph for
every $m$.

\begin{lem}\label{l9}
$B_m \sim B_n$ for $0 \le m \ne n \le p-1$.
\end{lem}
\begin{pf}
There exists a unique element $\hat{g} \in \kalapz \le P$  such that $\hat{g} (B_m) =B_n$. We also have a
unique element $\mathring{g} \in \mathring{\mathbb{Z}}_p \le P$
 with $\hat{g}^{\mathcal{B} }= \mathring{g}^{\mathcal{B} }$.
Since $\mathbb{Z}_p$ is cyclic and $\hat{x} \ne \mathring{x}$
we have $\hat{g} \ne \mathring{g}$.
Moreover, we may also assume
that $\hat{g}\restriction {B_m} \ne \mathring{g} \restriction {B_m}$ since
$\hat{g} \ne \mathring{g}$ and the induced subgraphs of $\Gamma$ by
$B_{m+c} \cup B_{n+c}$ are all isomorphic, where both $m+c$ and $n+c$ are taken modulo $p$.

Clearly, $\tilde{g} = \mathring{g} \hat{g}^{-1}$ is cycle of
length $p$ on $B_n$. The points of $V(\Gamma) \setminus \Lambda$
are contained in $P$-orbits of length $p$ so
 and $\tilde{g}$ fixes every point of the set
 $B_m \cup B_n \setminus \Lambda$ since $\tilde{g}^{\mathcal{B} }=id$.

Let $u \in B_m \setminus \Lambda$. It is enough to show that if $u$ is connected to some $v \in B_n$, then $u$ is
connected to every point of $B_n$. We will prove that $A$ is transitive on the following
pairs: $\left\{ \left(u,w \right) \mid w \in B_n \right\}$.

$A$ is transitive on $\left\{ \left(u,w \right) \mid w \in B_n \cap
  supp(\tilde{g}) \right\} =\left\{ \left(u,w \right) \mid w \in B_n \cap
  \Lambda \right\}$ since $\tilde{g}$ fixes $u$.
  Therefore we may assume that $v \in B_n \setminus \Lambda$ and we only have to find
  an element $a \in A$ such that $a(u) =u$ and $a(v) \in B_n \cap \Lambda$.

 The restriction of $\tilde{g}$ to $B_n$ is a cycle of length $p$
 which does not commute with $\mathring{r} \restriction B_n$,
 where $\mathring{r}$ is an involution of $\mathring{R}$.
 Since $\mathring{r}$ and $\mathring{g}$ commute we have
 that there is a $u' \in B_m$ such
that $\mathring{r} \hat{g}(u') \ne \hat{g} \mathring{r}
(u')$. Since the action of $\hat{R}$ is transitive on $B_m$
there exists
$\hat{r} \in \hat{R}$ such that $\hat{r}(u) =u'$.
Then \[ \left( \mathring{r} \hat{r} \right) \hat{g}  (u) =\mathring{r} \hat{g} \hat{r} (u) =\mathring{r} \hat{g}(u')
\ne \hat{g} \mathring{r} (u') = \hat{g} \left( \mathring{r} \hat{r}\right) (u)  \]
so there exists $a' \in A$ such that
\begin{equation}\label{eq1} a' \hat{g} (u) \ne \hat{g} a' (u) \mbox{.} \end{equation}

Let us assume that $v= \hat{g} (u)$. Then the inequality \ref{eq1} gives
  $a'(v) \ne \hat{g} a'(u) $.
Since $\hat{R} \restriction
  B_m$ is regular on $B_m$ there exists $\hat{s} \in \hat{R}$ such that
  $\hat{s}(u)=a'(u)$ and since $\hat{s} $ and $\hat{g}$
  commute we have $ \hat{s}(v) =  \hat{s} \hat{g} (u) = \hat{g} \hat{s} (u) = \hat{g} a' (u)$.
  Therefore $\hat{s} (v) \ne a'(v)$ and hence $\hat{s} ^{-1} a'$ fixes $u$ and $\hat{s}^{-1} a' (v) \ne v$ so we may assume that $v \ne \hat{g}(u)$.

If $p=7$, then $v \in B_n \cap \Lambda$.

If $p=5$, then there exists $\hat{t} \in \hat{R}$ such that $\hat{t}(u) \in B_m \setminus \Lambda =
B_m \setminus supp(\tilde{g}) $ while
$\hat{t}(v) \in B_n \cap \Lambda \subset supp(\tilde{g})$ since both $\hat{R}\restriction B_m$
 and $\hat{R}\restriction B_n$ are regular and $\gcd(8,5)=1$. The permutations
 $\hat{t}^{-1} \tilde{g}^l \hat{t}$
fix the point $u$  for every $0 \le l \le 4$ and $\hat{t}^{-1} \tilde{g}^{l_1} \hat{t}(y)
\ne \hat{t}^{-1} \tilde{g}^{l_2} \hat{t} (y)$ if $l_1 \not\equiv l_2$ $(mod \mbox{
  }p)$. At least one of the the four elements $\hat{t}^{-1} \tilde{g} \hat{t},
  \hat{t}^{-1} \tilde{g}^2 \hat{t},\hat{t}^{-1} \tilde{g}^3 \hat{t},\hat{t}^{-1} \tilde{g}^4 \hat{t}
$ of $A$ fixes $u$ and
  maps $v$ to an element of $B_n \cap fix(\tilde{g}) = B_n \cap \Lambda$ since $|B_n \setminus supp(\tilde{g})| =3$,
  finishing the proof of the fact that $B_m \sim B_n$ for $0 \le m \ne n
  \le 7$. \qed
\end{pf}

Every permutation of $V(\Gamma)$ which fixes setwise $B_m$ for every
$m$ is an automorphism of $\Gamma$ so there is an $a \in A$ such that
$\mathring{x}^a = \hat{x}$. Applying Proposition \ref{particio} we
get that there exists $\alpha \in A$ such that $\left(\hat{R} \times \hat{\mathbb{Z}}_p\right)
^{\alpha} = \mathring{R} \times \mathring{\mathbb{Z}}_p $,
finishing the proof of Proposition \ref{kicsi}.

\end{document}